\documentclass[10pt,letterpaper]{article}
\usepackage[top=0.85in,footskip=0.75in,marginparwidth=2in]{geometry}
\pdfoutput=1
\usepackage[utf8]{inputenc}

\usepackage{cite}

\usepackage{amsmath}

\usepackage{microtype}
\DisableLigatures[f]{encoding = *, family = * }

\raggedright
\usepackage{changepage}

\usepackage[aboveskip=1pt,labelfont=bf,labelsep=period,singlelinecheck=off]{caption}

\makeatletter
\renewcommand{\@biblabel}[1]{\quad#1.}
\makeatother

\usepackage{lastpage,fancyhdr,graphicx}
\usepackage{epstopdf}
\usepackage[symbol]{footmisc}
\fancyheadoffset[L]{2.25in}
\fancyfootoffset[L]{2.25in}

\usepackage{color}

\definecolor{Gray}{gray}{.25}

\usepackage{graphicx}

\usepackage{sidecap}

\usepackage{wrapfig}
\usepackage[pscoord]{eso-pic}
\usepackage[fulladjust]{marginnote}
\reversemarginpar

\begin{document}
\vspace*{0.35in}

% title goes here:

\begin{center}
{\Large
\textbf{High-order schemes for the Euler equations in cylindrical/spherical coordinates}
}
\newline
% authors go here:
\\
Sheng Wang\textsuperscript{1,}\footnote{Corresponding author, E-mail: sw767@cornell.edu},
Eric Johnsen\textsuperscript{2}
\\
\bigskip
1 Sibley School of Mechanical and Aerospace Engineering, Cornell  University, Ithaca, NY 14853, USA
\\
2 Department of Mechanical Engineering, University of Michigan, Ann Arbor, MI 48109, USA
\\

\bigskip

\end{center}
\justify
\section*{Abstract}
We consider implementations of high-order finite difference Weighted Essentially Non-Oscillatory (WENO) schemes for the
Euler equations in cylindrical and spherical coordinate systems with
radial dependence only. The main concern of this work lies in ensuring
both high-order accuracy and conservation. Three different spatial
discretizations are assessed: one that is shown to be high-order
accurate but not conservative, one conservative but not high-order
accurate, and a new approach that is both high-order accurate and
conservative.  For cylindrical and spherical coordinates, we present
convergence results for the advection equation and the Euler equations
with an acoustics problem; we then use the Sod shock tube and the
Sedov point-blast problems in cylindrical coordinates to verify our
analysis and implementations.

% the * after section prevents numbering
\section{Introduction}
In a variety of physical phenomena, the dominant dynamics occur in
spherical and cylindrical geometries.  Examples include astrophysics
(e.g., supernova collapse), nuclear explosions, inertial confinement
fusion (ICF) and cavitation-bubble dynamics.  A natural approach to
solving these problems is to write the governing equations in
cylindrical/spherical coordinates, which can then be solved
numerically using an appropriate discretization. Historically, the
first such numerical studies were conducted by Von Neumann and
Richtmyer \cite{Neumann1950} in the 1940s for nuclear explosions.  To
treat the discontinuities in a stable fashion, they explicitly
introduced artificial dissipation to the Euler equations.  While this
method correctly captures the position of shocks and satisfies the
Rankine-Hugoniot equations, flow features in the numerical solution,
in particular discontinuities, are smeared due to excessive
dissipation.

The collapse and explosion of cavitation bubbles, supernovae and ICF
capsules share similarities in that they are all, under ideal
circumstances, spherically symmetric flows that involve material
interfaces and shock waves. Such flows are rarely ideal, insofar as
they are prone to interfacial instabilities due to accelerations
(Rayleigh-Taylor \cite{Sharp1984}), shocks (Richtmyer-Meshkov
\cite{brouillette2002}), or geometry (Bell-Plesset
\cite{bell1951,plesset1956}).  When solving problems with large
three-dimensional perturbations, cylindrical/spherical coordinates may
not be advantageous. However, in certain problems such as
sonoluminescence \cite{Putterman1997}, the spherical symmetry
assumption is remarkably valid.  Modeling the bubble motion with
spherical symmetry can greatly reduce the computational cost.  As an
example, Akhatov, et al. \cite{Akhatov2001} used a first-order Godunov
scheme to simulate liquid flow outside of a single bubble whose radius
was given by the Rayleigh-Plesset equation. This approach assumes
spherical symmetry but does not solve the equations of motion inside
the bubble.

High-order accurate methods are becoming mainstream in computational
fluid dynamics \cite{wang2013}. However, implementation of such
methods in cylindrical/spherical geometries is not trivial.  Several
recent studies in cylindrical and spherical coordinates have focused
on the Lagrangian form of the equations
\cite{Maire2009,Vachal2014}. The Euler equations in cylindrical or
spherical geometry were studied by Maire \cite{Maire2009} using a
cell-centered Lagrangian scheme, which ensures conservation of
momentum and energy.  These equations were also considered by Omang et
al. \cite{Trulsen2006} using Smoothed Particle Hydrodynamics (SPH),
though SPH methods are generally not high-order accurate. On the other
hand, solving the equations in Eulerian form is not trivial,
especially when trying to ensure conservation and high-order
accuracy. Li \cite{Li2003} attempted to implement Eulerian finite
difference and finite volume weighted essentially non-oscillatory
(WENO) schemes \cite{Shu1996} on cylindrical and spherical grids, but
did not achieve satisfactory results in terms of accuracy and
conservation. Liu \emph{et al.} \cite{glaister1988} considered flux
difference-splitting methods for ducts with area variation.
\cite{liu1999} followed a similar formulation of the equations
employing a total variation diminishing method to simulate explosions
in air.  Johnsen \& Colonius \cite{johnsen2008,johnsen2009} used
cylindrical coordinates with azimuthal symmetry to simulate the
collapse of an initially spherical gas bubble in shock-wave
lithotripsy by solving the Euler equations inside and outside the
bubble using WENO.  De Santis \cite{desantis2014} showed equivalence
between their Lagrangian finite element and finite volume schemes in
cylindrical coordinates.  Xing and Shu
\cite{Xing2005,Xing2006_1,Xing2006_2} performed extensive studies of
hyperbolic systems with source terms, which are relevant as the
equations in cylindrical/spherical coordinates can be written with
geometrical source terms.  Although one of their test cases involved
radial flow in a nozzle using the quasi one-dimensional nozzle flow
equations, they did not consider the general gas dynamics
equations. Thus, at this time, a systematic study of the Euler
equations in cylindrical and spherical coordinates, with respect to
order of accuracy and conservation, has yet to be conducted.

In this paper, we investigate three different spatial discretizations
in cylindrical/spherical coordinates with radial dependence only using
finite difference WENO schemes for the Euler equations. In particular, we propose a new approach that is both
high-order accurate and conservative. Here, we are concerned with the
interior scheme; appropriate boundary approaches will be investigated
in a later study. The governing equations are stated in Section
\ref{sec:numerical framework} and the spatial discretizations are
presented in Section \ref{sec:Numerical method}. In Section
\ref{sec:numerical results}, we test the different discretizations on
smooth problems (scalar advection equation, acoustics problem for the
Euler equations) for convergence, and with shock-dominated problems
(Sod shock tube and Sedov point blast problems) for conservation. The
last section summarizes the present work and provides a future
outlook.

\section{Governing equations}
  \label{sec:numerical framework}
The differential for the Euler equations in cylindrical/spherical coordinates
  with radial dependence only:
%  \begin{subequations}\label{eq:compressible_euler} \begin{align}
 % \left ( \rho \right )_{t}+\frac{1}{r^{\alpha }}\left ( r^{\alpha
 % }\rho u \right )_{r}&=0, \\ \left ( \rho u\right
 % )_{t}+\frac{1}{r^{\alpha }}\left ( r^{\alpha }\rho u^{2}
 % \right)_{r}+p_{r}&=0, \\ \left ( E \right )_{t}+\frac{1}{r^{\alpha
 % }}\left[ r^{\alpha }\left ( E+p \right)u \right]_{r}&=0,
 % \end{align} \end{subequations} 	
 
\begin{align}  \label{eq:compressible_euler} 
 \frac{\partial \mathbf{U}}{\partial t}+ \frac{1}{r^{\alpha}}\frac{\partial \left ( r^{\alpha}\mathbf{F}(\mathbf{U}) \right )}{\partial r}=\mathbf{S}(\mathbf{U})
  \end{align} 
  where $\mathbf{U}=( \rho, m=\rho u, E )^{T}$, $\mathbf{F}(\mathbf{U}) = \left ( m,    \frac{m^2}{\rho}+p,   \frac{m}{\rho}(E+p) \right )^{T}$, and 
  $\mathbf{S}(\mathbf{U}) = (0,p,0)^T$. 
  $\rho$ is the density, u is velocity in radial direction,  
  $t$ is time, $r$ is the radial coordinate, $p$ is the
  pressure, $E$ is the total energy per unit
  volume, and $\alpha$ is a geometrical parameter, which is 0, 1, or 2
  for Cartesian, cylindrical, or spherical coordinates, respectively.
  Subscripts denote derivatives.  Diffusion effects are neglected. For
  an ideal gas, the equation of state to close this system can be
  written:
  \begin{equation} p=\left
  ( \gamma -1 \right )\varepsilon, \label{eq: Equation of state}
  \end{equation} 
  where $\varepsilon =E-\rho u^{2}/2$ is the internal energy per unit
  volume, and $\gamma$ is the specific heats ratio.  Other equations
  of state can be used, e.g., a stiffened equation for liquids and solids.

\section{Numerical method} \label{sec:Numerical method}

We describe three discretizations of the Euler
Eqs.\eqref{eq:compressible_euler} in cylindrical/spherical
 coordinates that differ based on the treatment of the convective
 terms.  While the discretized form of the Euler equations in
 Cartesian coordinates is generally designed to conserve mass,
 momentum and energy, the conservation condition does not necessarily
 hold in cylindrical or spherical coordinates, depending on the
 numerical treatment of the equations.  The criterion we use to
 determine discrete conservation is as follows:
  \begin{equation} \label {eq: conservation}
    \int_{\Omega}\phi(t,r)dv=\int_{\Omega}\phi(0,r)dv,
  \end{equation}
  where $\Omega$ is a domain (possibly a computational cell).  Here,
  $\phi$ represents a conserved variable, e.g., density, momentum per
  unit volume or energy per unit volume. This equation means that the
  total mass, momentum, and energy are constant in time provided there
  is no flux of these quantities through the boundaries of the domain,
  which is the case for the problems of interest here.  A different
  approach to defining conservation for hyperbolic laws is the exact
  C-type property \cite{Xing2005}, which implies that the system
  admits a stationary solution in which nonzero flux gradients are
  exactly balanced by the source terms in the steady-state case.  Xing
  and Shu \cite{Xing2006_1,Xing2006_2} applied WENO in systems of
  conservation laws with source terms and considered radial flow in a
  nozzle using the quasi one-dimensional nozzle flow equations. In our
  work, we focus on the Euler equations.
  
  We consider finite difference(FD) WENO
  schemes. we give  brief description of a fifth-order finite difference WENO scheme is given in Appendix.  For finite difference WENO, given the cell-centered values, the fluxes are first split and then interpolated to compute the numerical flux. In order to give a clear image of implementation, we write the solution procedure right after describing the spatial discretization for each method.  

\subsection*{Method One}
 The first spatial discretization, labelled \emph{Method One} here,
 can be found in Chapter 1.6 of Toro \cite{toro1999}. Expand the convective term and move the part without spatial derivative to right hand of Eq.~\ref{eq:compressible_euler} to obtain the differential equations 
 
 \begin{align}  \label{eq:differential_FD1} 
 \frac{\partial \mathbf{U}}{\partial t}+\frac{\partial \mathbf{F}(\mathbf{U})}{\partial r}=\mathbf{S}(\mathbf{U})- \frac{\alpha}{r}\mathbf{F}(\mathbf{U})
  \end{align} 
 The notation is consistent with the notation in Eq.~\ref{eq:compressible_euler}. 
 
 For the convenience of programming, we also write the mass, momentum and energy in semi-discrete form:
 \begin{subequations} \label{eq: compressible_euler_FD1}
  \begin{align}
   \frac{\mathrm{d}\rho_{i} }{\mathrm{d} t}=&- \frac{ (\rho
  u)_{i+1/2}- (\rho u)_{i-1/2}}{\Delta r_i} - \frac{\alpha}{r_{i}}\left
  ( \rho u \right )_{i}, \\ \frac{\mathrm{d}(\rho u)_{i} }{\mathrm{d}
  t}=& - \frac{ (\rho u+p)_{i+1/2}- (\rho u+p)_{i-1/2}}{\Delta r_i} -
  \frac{\alpha}{r_{i}}\left ( \rho u ^{2}\right )_{i}, \\
  \frac{\mathrm{d} E_{i} }{\mathrm{d} t}=& - \frac{ [(E+p) u]_{i+1/2}-
  [(E+p) u]_{i-1/2}}{\Delta r_i} - \frac{\alpha}{r_{i} }\left[ u( E+p)
  \right]_{i},
  \end{align}
  \end{subequations}  
  where $\Delta r_i$ is the linear radial cell width. For FD WENO, the variables to be
  evolved in time are the cell-centered values, i.e., the values at
  $r_{i}$ in cell $[r_{i-1/2},r_{i+1/2}]$.

  With this approach derived from the differential form of the
  equations rather than the integral, physical variables expected to
  be conserved are not necessarily conserved numerically because there
  is no strict constraint between the conservative fluxes and the
  geometrical source terms.  However, high-order accuracy may be achieved with
  this method.\\
  
  This part summarizes the solution procedure for Method One used in our code:\\
1. Initialize the primitive variables, $\rho$, $u$, $p$, and $E$ \\
2. Using local Lax-Fredrich to split the flux, obtain $\rho u_i^{\pm}$, $(\rho u+p)_i^{\pm}$ , $[(E+p) u]_{i}^{\pm}$, and local speed of wind, $\lambda_i$. The plus sign means the flux moves toward right, the minus sign means the flux moves toward left. The convention is used in all the flux calculation part in this paper.\\
3. Using the local characteristic decomposition and finite difference WENO to approximate the flux, obtain $\rho u_{i+1/2}^{\pm}$ ,  $(\rho u+p)_{i+1/2}^{\pm}$ , and $[(E+p) u]_{i+1/2}^{\pm}$.\\
4. Calculate the residual. The source terms in Method one are collocated with primitive variable, can be directly added to the residual. Note that the source terms are updated in each sub-step.\\
5. March in time.\\

\subsection*{Method Two}

 The second discretization, Method Two, is based on the integral form
 of the equations and was used by several authors
 \cite{glaister1988,liu1999,johnsen2008}. The mass, momentum and
 energy equations are written in semi-discrete form:
  \begin{subequations}\label{eq: compressible_euler_mass_FD2}
  \begin{align}
  \frac{\mathrm{d}\rho_{i} }{\mathrm{d} t}=&-\frac{r_{i+1/2}^{\alpha
    }(\rho u)_{i+1/2}-r_{i-1/2}^{\alpha}(\rho u)_{i-1/2}}{\Delta V_i}, \\ 
    \frac{\mathrm{d}(\rho u)_{i} }{\mathrm{d} 
      t}=&-\frac{r_{i+1/2}^{\alpha }(\rho
        u^{2}+p)_{i+1/2}-r_{i-1/2}^{\alpha}(\rho u^{2}+p)_{i-1/2}}{\Delta V_i}
	  +S(r_i), \\ 
	     \frac{\mathrm{d} E_{i} }{\mathrm{d} t}=&-\frac{r_{i+1/2}^{\alpha
	         }[(E+p)u]_{i+1/2}-r_{i-1/2}^{\alpha}[(E+p)u]_{i-1/2}}{\Delta V_i},
   \end{align}
   \end{subequations}
   where $\Delta V_i=\frac{1}{1+\alpha
   }(r_{i+1/2}^{\alpha+1}-r_{i-1/2}^{\alpha+1 })$ and $S(r_i)$ is the
   source term in the momentum equation, which can be expressed as:
 \begin{equation}  
 S(r_i) = \frac{r_{i+1/2}^{\alpha
 }p_{i+1/2}-r_{i-1/2}^{\alpha }p_{i-1/2}}{\Delta
 V_i}-\frac{p_{i+1/2}-p_{i-1/2}}{\Delta r} .
 \end{equation} 
 Depending on the reconstruction procedure, the first term may cancel
 the corresponding term in the momentum equation.  Again, the
 variables to be evolved in time are the cell-centered values for FD WENO.

 With this approach, the relevant physical variables are expected to
 be conserved. However, the order accuracy is ultimately second
 order. This latter point can be readily understood by subtracting
 Method Two from Method One:
\begin{equation}
\frac{f_{j+1/2 }r_{j+1/2 }^{\alpha}-f_{j-1/2 }r_{j-1/2}^{\alpha}}{\frac{1}{\alpha +1}(r_{j+1/2}^{\alpha+1}-r_{j-1/2}^{\alpha+1})}- \left(\frac{f_{j+1/2}-f_{j-1/2}}{\Delta r}+\frac{\alpha f_{j}}{r_{j}} \right)=O(\Delta r^{2}).
\end{equation}  
   
  This part summarizes the solution procedure for Method Two used in our code: \\
1. Initialize the primitive variables, $\rho$, $u$, $p$, and $E$ \\
2. Using local Lax-Fredrich to split the flux, obtain $\rho u_i^{\pm}$, $(\rho u+p)_i^{\pm}$ , $[(E+p) u]_{i}^{\pm}$, and local speed of wind, $\lambda_i$. \\
3. Using the local characteristic decomposition and finite difference WENO to approximate the flux, obtain $\rho u_{i\pm1/2}^{\pm}$ ,  $(\rho u+p)_{i\pm1/2}^{\pm}$ , and $[(E+p) u]_{i\pm1/2}^{\pm}$.\\
4. Calculate the residual. The source terms in Method Two are located at the cell surface. It is not straight forward to calculate it form $\rho u_{i\pm1/2}^{\pm}$ ,  $(\rho u+p)_{i\pm1/2}^{\pm}$ , and $[(E+p) u]_{i\pm1/2}^{\pm}$. Our solution is to approximate $p_{i \pm 1/2}^{\pm}$ using  the same nonlinear weights.  The source terms are updated in each sub-step.\\
5. March in time.\\

\subsection*{Method Three}
The third spatial discretization, Method Three, is inspired by the
solution to acoustics problems in spherical coordinates. This approach is also used by Toro \cite{toro1999} and  Zhang\cite{Zhang2011}. Multiplying Eqs.~\eqref{eq:compressible_euler} by $r^{\alpha }$,
 \begin{align}  \label{eq:differential_FD3} 
\frac{\partial (r^{\alpha}\mathbf{U})}{\partial t}+\frac{\partial \left ( r^{\alpha}\mathbf{F}(\mathbf{U}) \right )}{\partial r}=r^{\alpha}\mathbf{S}(\mathbf{U})
  \end{align}

For the convenience of programming, the mass, momentum, and energy equations are written in semi-discrete form:
  \begin{subequations}\label{eq:compressible_euler_mass_FD3}
  \begin{align}
\frac{\mathrm{d}(r^{\alpha }\rho)_{i} }{\mathrm{d}
t}=&-\frac{(r^{\alpha }\rho u)_{i+1/2}-(r^{\alpha} \rho
u)_{i-1/2}}{\Delta r_i}, \\ 
\frac{\mathrm{d}(r^{\alpha }\rho u)_{i} }{\mathrm{d} t}=&-
\frac{[r^{\alpha }(\rho u^{2}+p)]_{i+1/2}-[r^{\alpha }(\rho
u^{2}+p)]_{i-1/2}}{\Delta r_i}+\alpha(pr^{\alpha -1})_{i}, \\ 
\frac{\mathrm{d}(r^{\alpha }E)_{i} }{\mathrm{d} t}=& -\frac{[r^{\alpha
}(E+p)u]_{i+1/2}-[r^{\alpha }(E+p)u]_{i-1/2}}{\Delta r_i}.  \\ 
  \notag \end{align}
\end{subequations}  
 For FD WENO, the cell-centered values are considered. 
This approach strictly follows the integral form of the Euler equations in cylindrical or spherical coordinates and
satisfies the C-type property for hyperbolic equations with source
terms \cite{Xing2005}. Thus, it is expected to be both conservative
and high-order accurate.

   \begin{figure}[t]% order of placement preference: here, top, bottom
  \centering
  \includegraphics[width=0.5\textwidth]{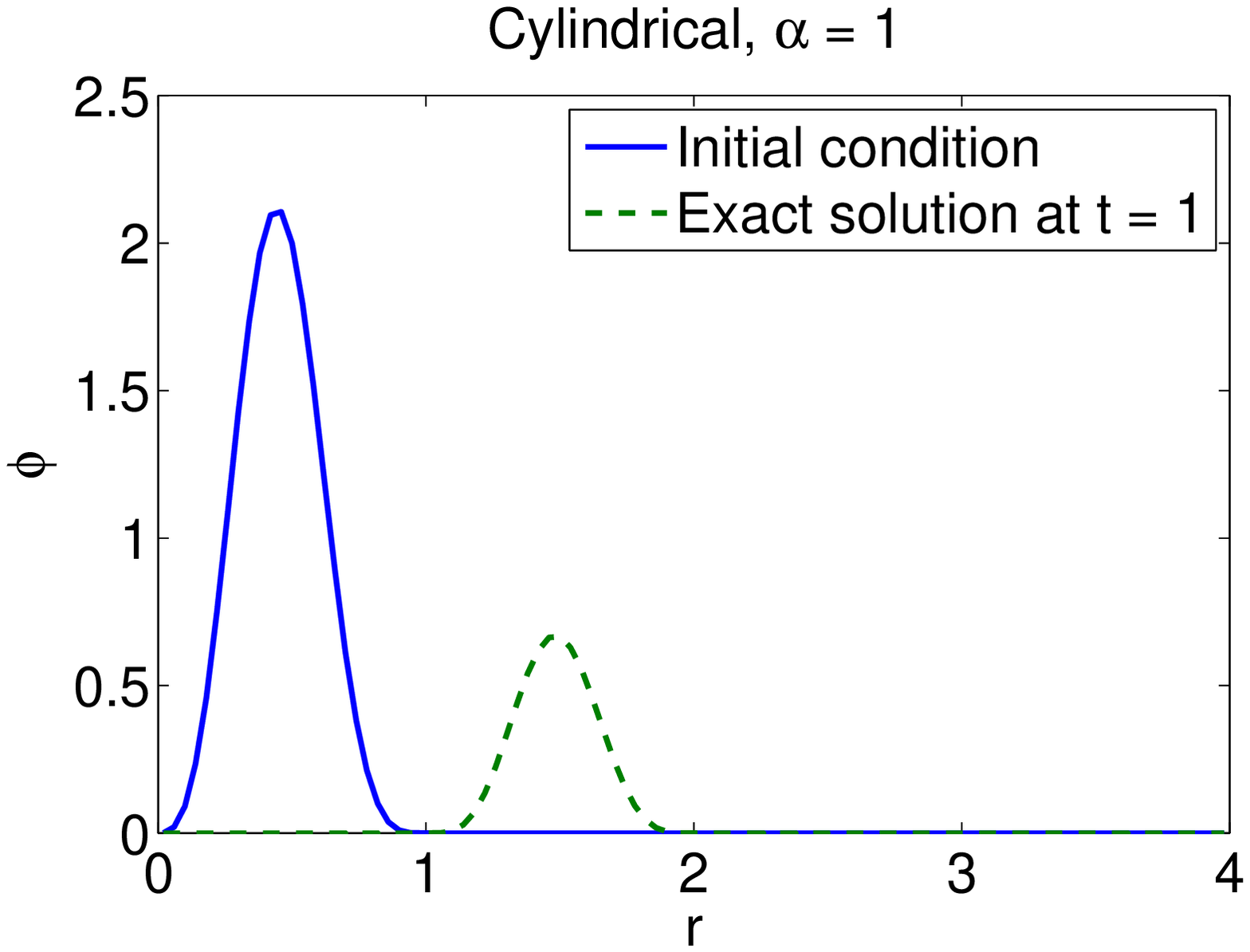}
  \includegraphics[width=0.48\textwidth]{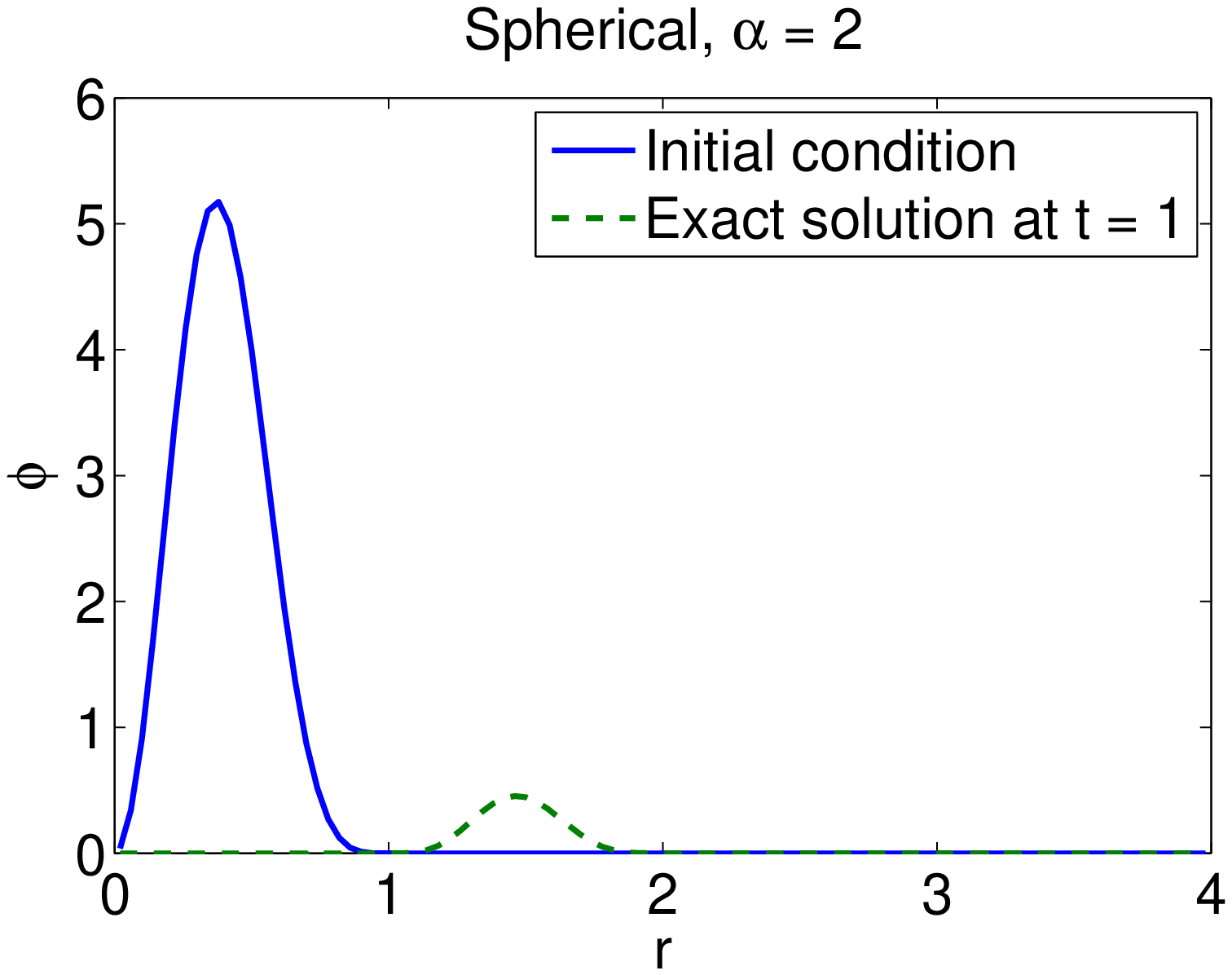}
  \caption{Initial conditions and exact solution after $t =1$ for the advection equation. }
  \label{fig:Initial condition for advection equation}
\end{figure}

  This part summarizes the solution procedure for Method Two used in our code: \\
1. Initialize the primitive variables, $\rho$, $u$, $p$, and $E$, then calculate $r^\alpha\rho$, $r^\alpha u$, $r^\alpha p$, and $r^\alpha E$ for each cell. \\
2. Using local Lax-Fredrich to split the flux, obtain $(r^\alpha\rho u)_i^{\pm}$, $[r^\alpha(\rho u^2+p)]_i^{\pm}$ , $[r^\alpha(E+p) u]_{i}^{\pm}$, and local speed of wind, $\lambda_i$. \\
3. Using the local characteristic decomposition and finite difference WENO to approximate the flux, obtain $(r^\alpha\rho u)_{i\pm 1/2}^{\pm}$, $[r^\alpha(\rho u^2+p)]_{i\pm 1/2}^{\pm}$ , $[r^\alpha(E+p) u]_{i\pm 1/2}^{\pm}$.\\
4. Calculate the residual. The source terms in Method Three are collocated with primitive variable, can be directly added to the residual. The source terms are updated in each sub-step.\\
5. March in time.\\

\section{Numerical results}
\label{sec:numerical results}
In this section, we apply the three discretizations introduced in the
previous section to four test cases using fifth-order WENO in
characteristic space with Local Lax-Friedrichs upwinding, and
fourth-order explicit Runge-Kutta with a Courant number of $0.5$ for
time marching. First, we use two smooth problems (scalar advection and
acoustics for the Euler equations) to demonstrate the convergence
rates of each method for the interior solution, with no regards
for boundary schemes. Next, we test conservation with two
shock-dominated problems (Sod shock tube and Sedov point blast
problems).

\subsection{Scalar advection problem}
Before considering nonlinear systems, the scalar advection equation is 
investigated. The advection equation in cylindrical and spherical coordinates with
symmetry is  written: 
 \begin{align} \left ( \phi \right
 )_{t}+\frac{c_{0}}{r^{\alpha }}\left ( r^{\alpha }\phi \right )_{r}=0,
 \label{eq: advection_equation} \end{align} 
where $\phi$ is a scalar field, $c_{0}$ is the (constant and known) wave
speed. Here, $c_0 =1$.
The initial conditions are	 
\begin{align}
   \phi(r,0)= \begin{cases}
\frac{\sin^{4}(\pi r) }{r^{\alpha }}, & \text{ if } 0\leq r\leq 1, \\ 
 0, & \text{ if } r> 1.
 \end{cases}
 \end{align}
 
 For this problem, the exact solution at time $t=1$ is
   \begin{align}
\phi(r,t)= \begin{cases}
\frac{\sin^{4}(\pi (r-c_{0}t)) }{r^{\alpha }}, & \text{ if }c_{0}t\leq r\leq c_{0}t+1, \\ 
 0, & \text{ otherwise}.
 \end{cases}
 \end{align}
The initial conditions and exact solution at $t=1$ are shown in
Fig.~\ref{fig:Initial condition for advection equation}.  Nearly
identical set-ups are used for the cylindrical and spherical cases,
the only difference being the geometrical parameter: $\alpha=1$ for
the cylindrical case, and $\alpha=2$ for the spherical.

\begin{figure}[t]% order of placement preference: here, top, bottom
\centering
\includegraphics[width=0.8\textwidth]{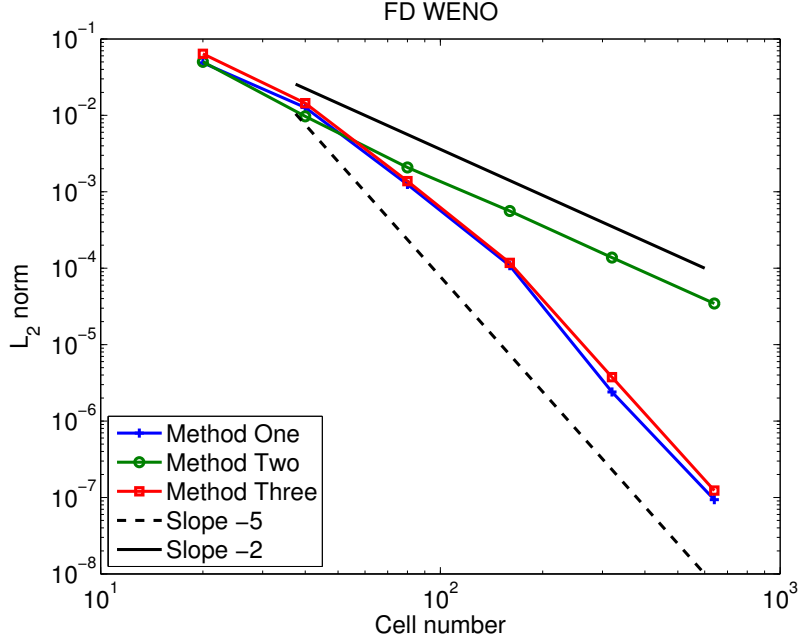}
\caption{$L_{2}$ error norm for all three discretizations for the scalar advection problem.}
\label{fig:advection_second_norm}
\end{figure}

The goal is to determine the convergence rates of each method
independently of boundary schemes. The problem set-up is specifically
chosen to prevent any boundary effects.  Here, we show convergence
results only for cylindrical coordinates, as the convergence rate is
similar for the spherical case.  Grids with $N= 21, 41, 81, 161, 321,
641$ are considered with constant $\Delta r$, and the exact solution
is used to evaluate the error of each
solution. Fig.~\ref{fig:advection_second_norm} shows the $L_2$ error
norm to verify the order of accuracy.  Methods One and Three both
achieve close to fifth-order accuracy, while Method Two is only
second-order accurate, as expected from the discussion in the previous
section.

\begin{figure}[t]% order of placement preference: here, top, bottom
\centering
\includegraphics[width=0.8\textwidth]{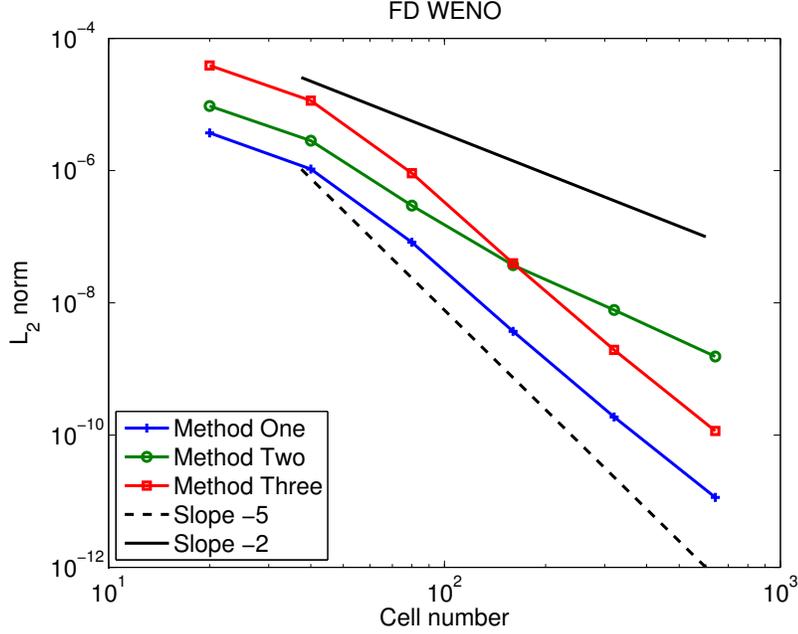}
\caption{$L_{2}$ error norm in density for each discretization for the
  acoustics problem. }
\label{fig:Euler_second_norm}
\end{figure}  
 
\subsection{Euler equations: acoustics problem  }

 A smooth problem is used to verify convergence with the Euler
 equations. The acoustics problem from Johnsen \&
 Colonius \cite{Eric2006} is adapted to spherical coordinates, with
 initial conditions:
   \begin{subequations}
 \begin{align}
\rho (r,0)=&1+\varepsilon f(r), \\ 
u(r,0)=&0, \\
 p(r,0)=&1/\gamma +  \varepsilon f(r), 
 \end{align}
 \end{subequations}
 with perturbation
 \begin{equation}
 f(r)= \begin{cases}
\frac{\sin^{4}(\pi r)}{r},& \text{ if } 0.4\leq r\leq 0.6, \\ 
 0, &  \mbox{otherwise}.
 \end{cases}
 \end{equation}
  For a sufficiently small $\varepsilon$ (here $10^{-4}$), the
  solution remains very smooth.  In this problem, the initial
  perturbation splits into two acoustic waves traveling in opposite
  directions. To prevent the singularity at the origin and boundary
  effects, the final time is set such that the wave has not yet
  reached the origin. Again, grids with $N=$ 21, 41, 81, 161, 321 and
  641 are used with constant $\Delta r$. Although an exact solution to
  order $\varepsilon^2$ is known, the solution on the finest grid is
  used as the reference to evaluate the error.

  Fig.~\ref{fig:Euler_second_norm} shows the $L_2$ error in density
  for this problem. The results show that Methods One and Three remain
  high-order and in fact fall on top of each other, although the rate
  now is closer to fourth order. Again, for Method Two, the rate is
  second order, as expected.

  \subsection{Sod shock tube}
    
  We consider the Sod shock tube problem \cite{sod1978} in cylindrical
  coordinate with azimuthal symmetry. The initial conditions are:
  \begin{equation}
  \left(  \begin{array}{c}
  \rho \\
  u \\
  p \end{array}\right)_L = 
  \left(  \begin{array}{c}
  1 \\
  0 \\
  1 \end{array}\right), \qquad
  \left(  \begin{array}{c}
  \rho \\
  u \\
  p \end{array}\right)_R = 
  \left(  \begin{array}{c}
  0.125 \\
  0 \\
  0.1 \end{array}\right).
  \end{equation}
  The domain size is 1 and 100 equally spaced grid points are
  used. The location of the ``diaphragm'' separating the left and
  right states is $r=0.5$.  Since no wave reaches the boundaries over
  the duration of the simulation (final time: $t_{f}=0.2$), the
  boundary scheme is irrelevant.

 \begin{figure}[t]% order of placement preference: here, top, bottom
 \centering
 \includegraphics[width=1\textwidth]{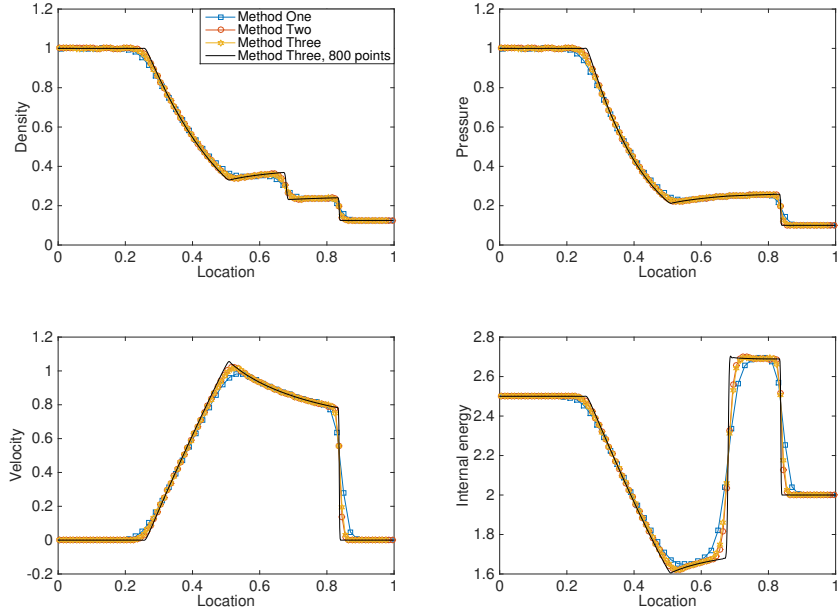}
 \caption{Profiles at $t=0.2$ for the Sod problem with 100 points.} 
 \label{fig:Comparing Sod}
 \end{figure} 

\begin{figure}[t]% order of placement preference: here, top, bottom
\includegraphics[width=0.46\textwidth]{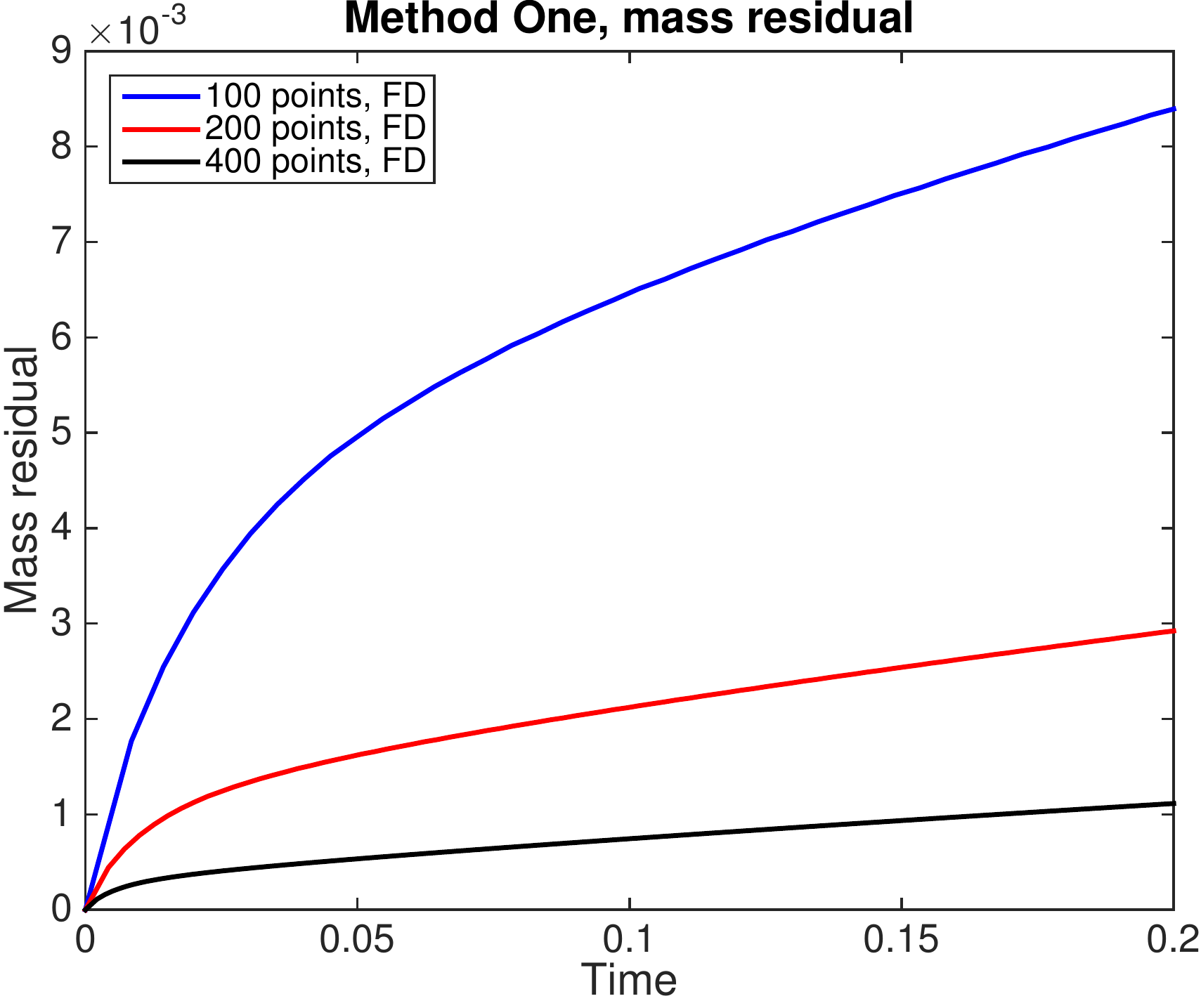}\includegraphics[width=0.47\textwidth]{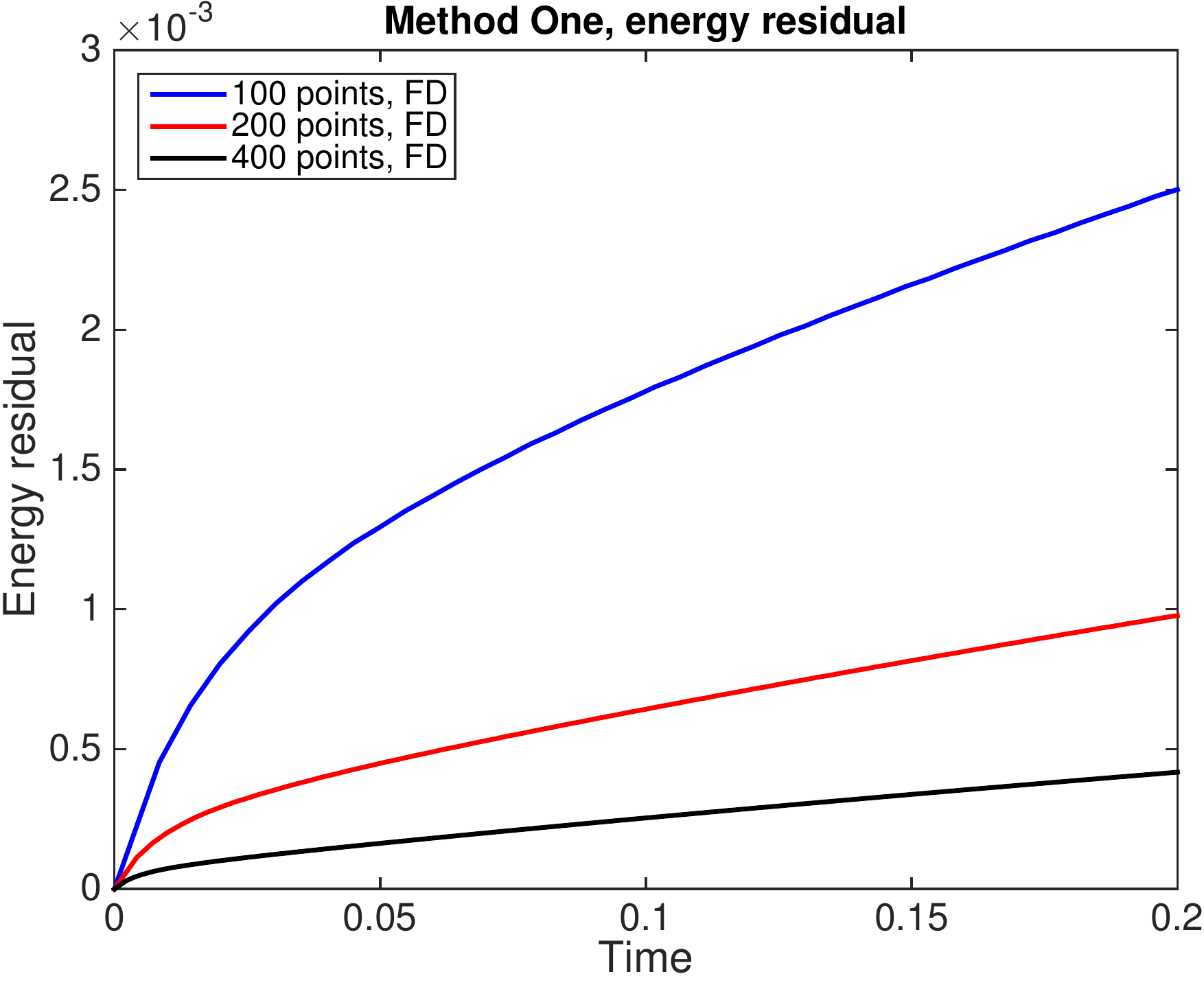}
\includegraphics[width=0.46\textwidth ]{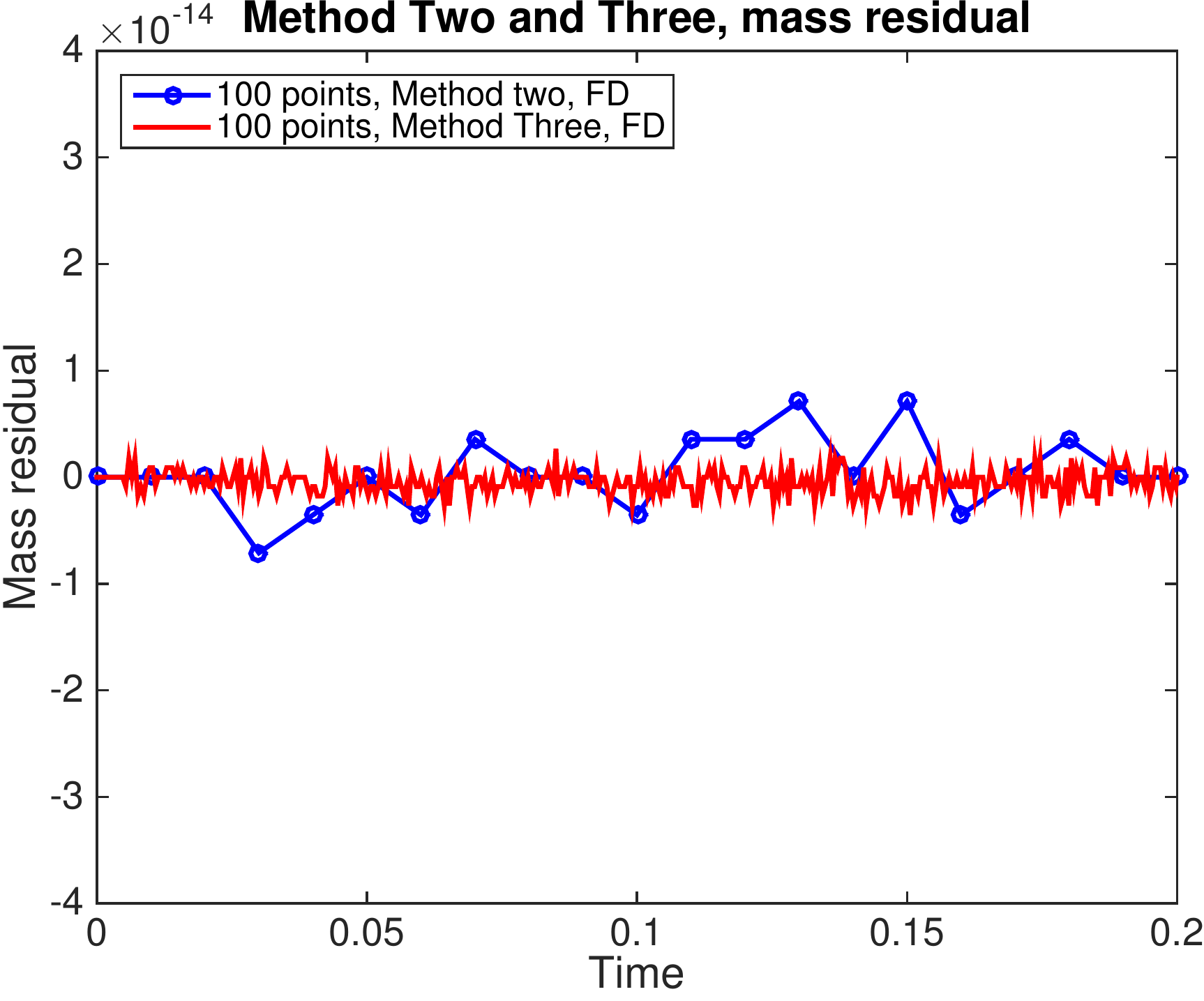}    \includegraphics[width=0.46\textwidth]{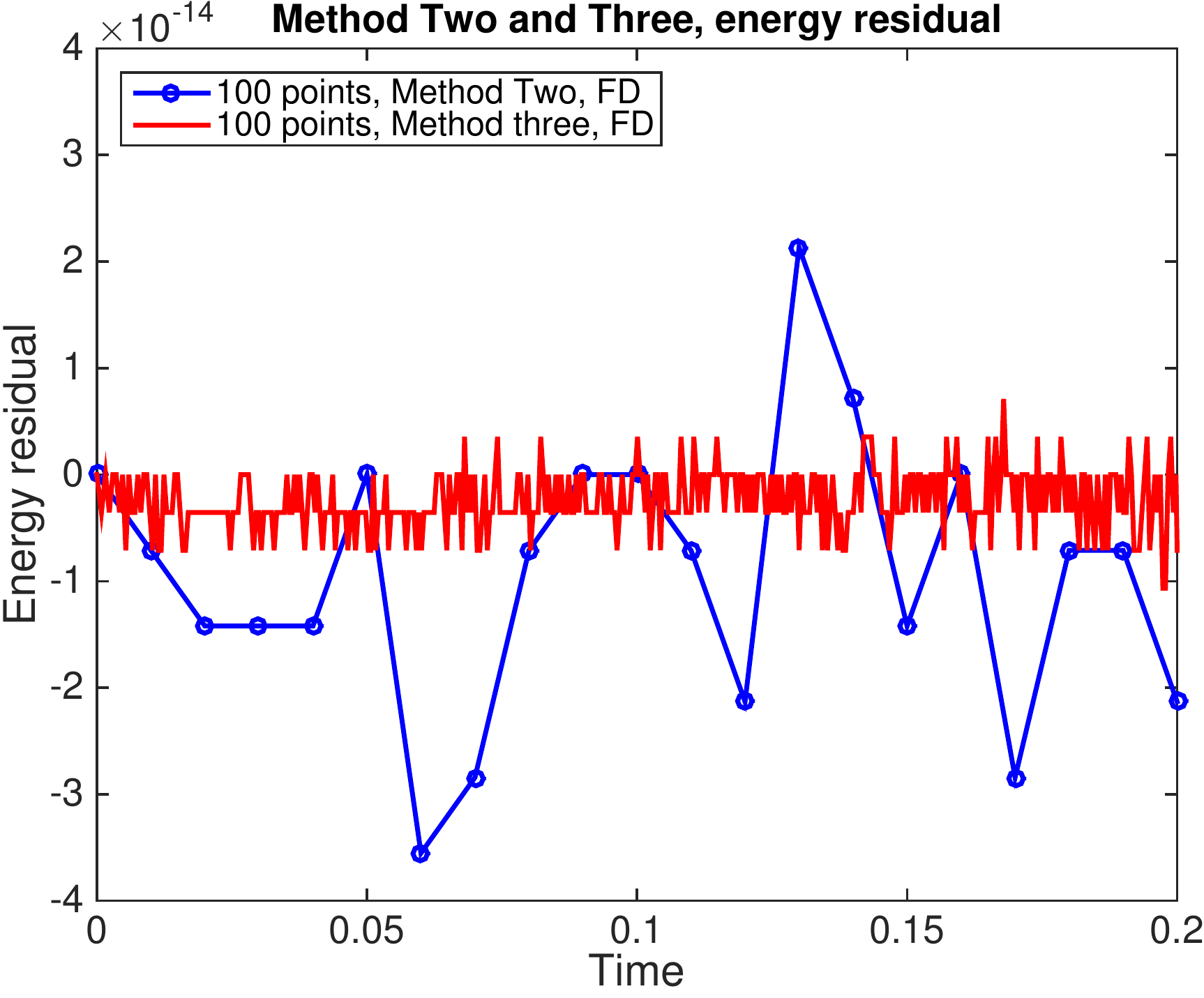}
\caption{Mass and energy residuals vs. time for the Sod problem.}
\label{fig:Residual of mass Sod}
\end{figure}

  Fig.~\ref{fig:Comparing Sod} shows density, velocity, pressure and
  internal energy profiles for this problem at the final time. Method Three on a gird of 800 points is  used as a reference solution. 
  On this grid, all three methods produce similar profiles. However, the
  residuals of the total mass and energy,
  \begin{equation}
  \Delta = \int_{I}\phi(t,r)dv-\int_{I}\phi(0,r)dv,
  \end{equation}
  yield different results, as observed in Fig.~\ref{fig:Residual of
    mass Sod}. For FD WENO, a high-order Gauss quadrature is employed to
  integrate the total mass and total energy from the cell-centered
  values.  As shown in Fig.~\ref{fig:Residual of mass Sod}, while
  Methods Two and Three are conservative to round-off level, Method
  One is not discretely conservative, as expected. In
  Fig.~\ref{fig:Comparing Sod}, differences in shock position due to
  lack of conservation are not clear.  The total momentum is zero
  based on the azimuthally symmetric setup.

\subsection{Sedov point-blast}

 Finally, we consider the Sedov point-blast problem in spherical
 coordinates.  Following the set-up of Fryxell et
 al. \cite{Fryxell2000} the initial conditions are
 \begin{equation}
  \left(  \begin{array}{c}
  \rho \\
  u \\
  p \end{array}\right) = 
  \left(  \begin{array}{c}
  1 \\
  0 \\
  1 \end{array}\right),
 \end{equation}
 except for a few computational cells around the origin, whose pressure is
  \begin{align} p_{0}^{'}=\frac{3(\gamma -1)\varepsilon}{(\alpha +2)\pi \delta r^{\alpha+1}}.
   \end{align}
   Here, $ \varepsilon =0.44$ is the dimensionless energy per unit
   volume. $\gamma$ is the specific heats ratioand and $\alpha$  a geometrical parameter, which is consistent with the value in Section \ref{sec:numerical framework}.  The domain size is 1, and $N=100$ with uniform spacing. We
   choose a constant $\delta r$ to be three times as large as the cell
   size for $N=100$.  Reflecting boundary conditions
   \cite{mohseni2000} are used along the centerline; since the shock
   does not leave the domain, the outflow boundary scheme is irrelevant.
   Due to the reflecting boundary condition at the center, the high
   pressure region is made up of 6 cells, i.e., 3 ghost cells and 3
   cells in the interior. The solution is plotted at $t=1$.

       \begin{figure}[t]% order of placement preference: here, top, bottom
   \centering
   \includegraphics[width=1\textwidth]{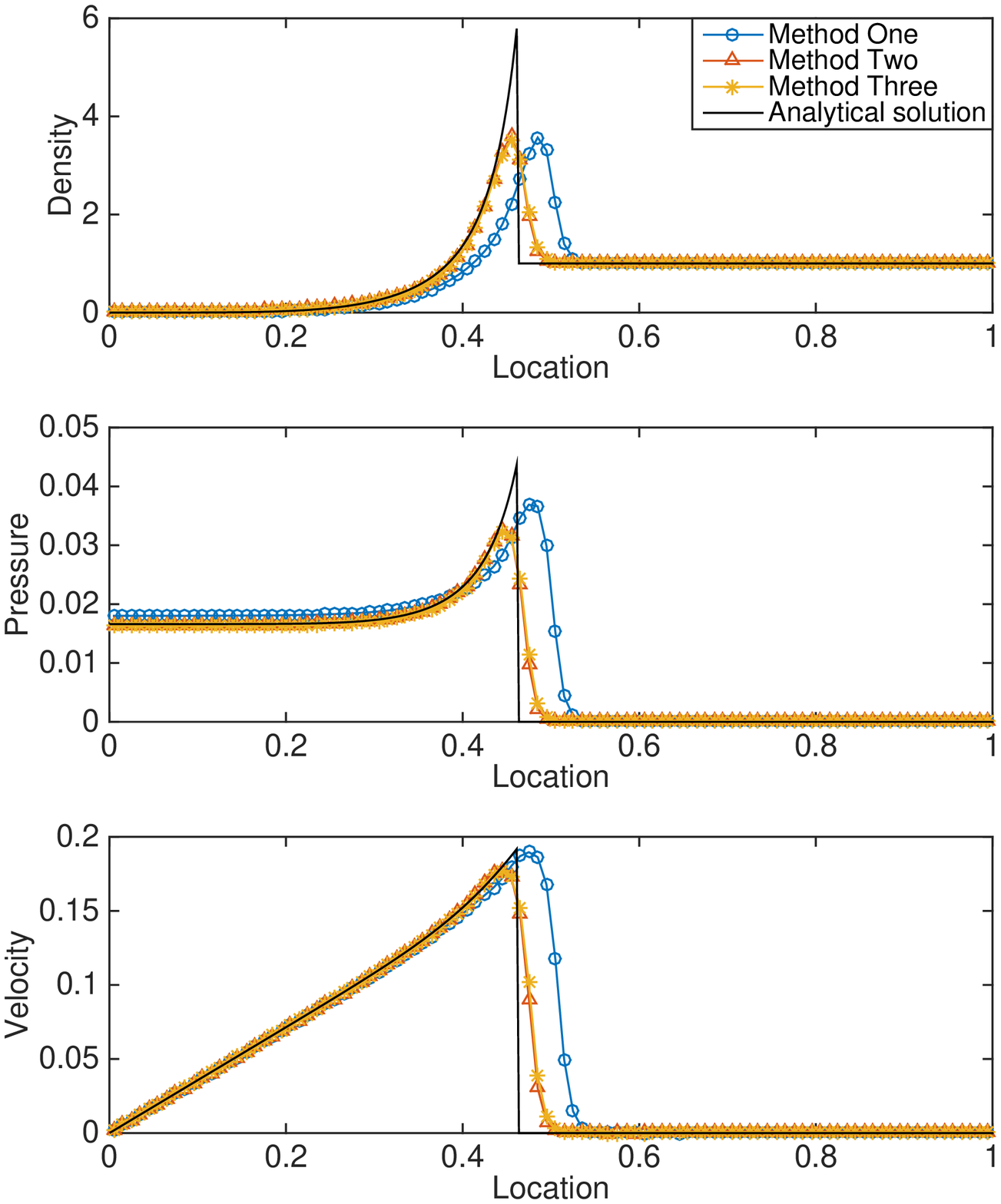}
   \caption{Profiles at $t=1$ for the Sedov problem with 100 points, FD WENO.}
   \label{fig:FD_Comparing Sedov}
   \end{figure}
   
\begin{figure}[ht]% order of placement preference: here, top, bottom
\centering
\includegraphics[width=0.47\textwidth]{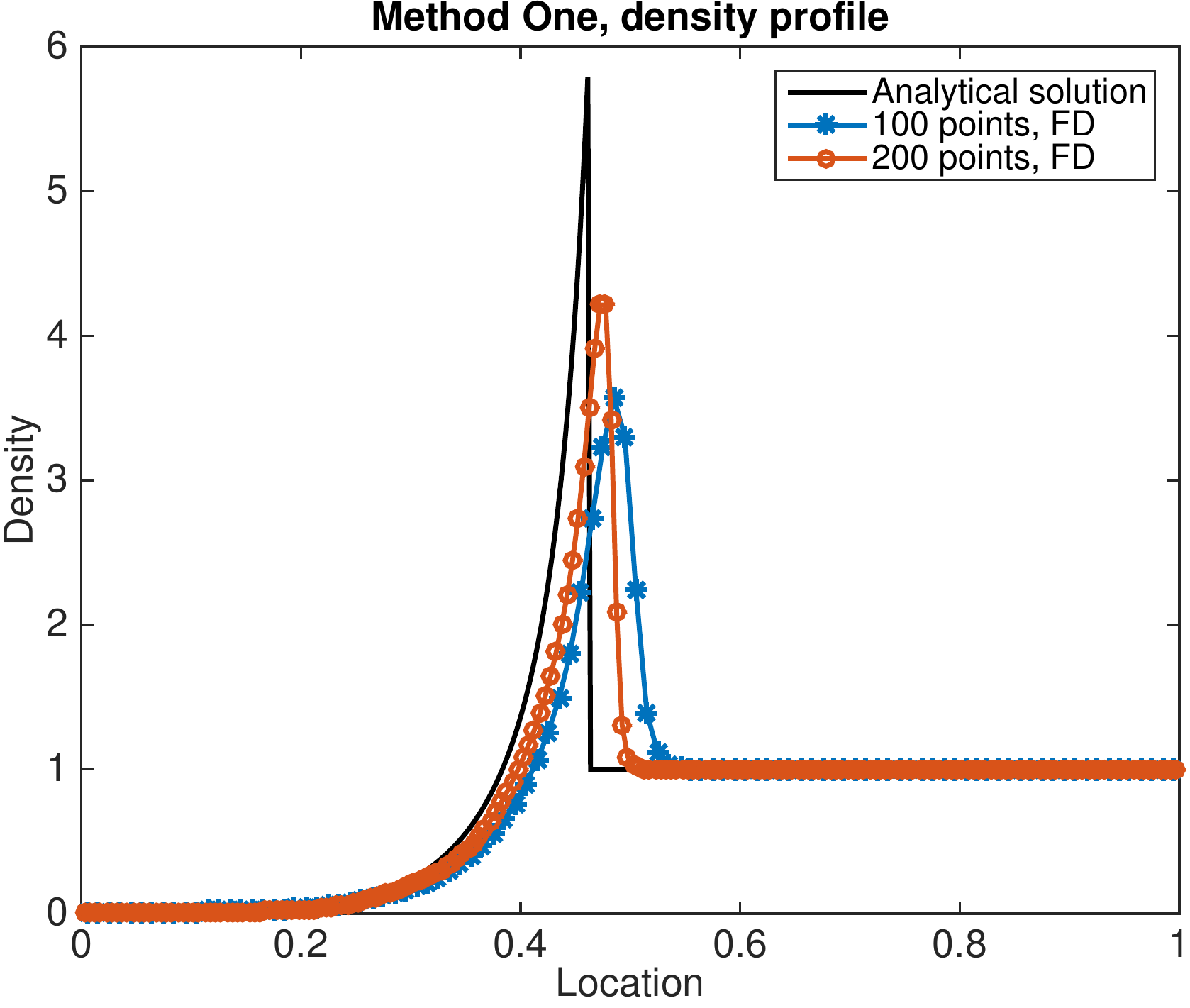} \includegraphics[width=0.46\textwidth]{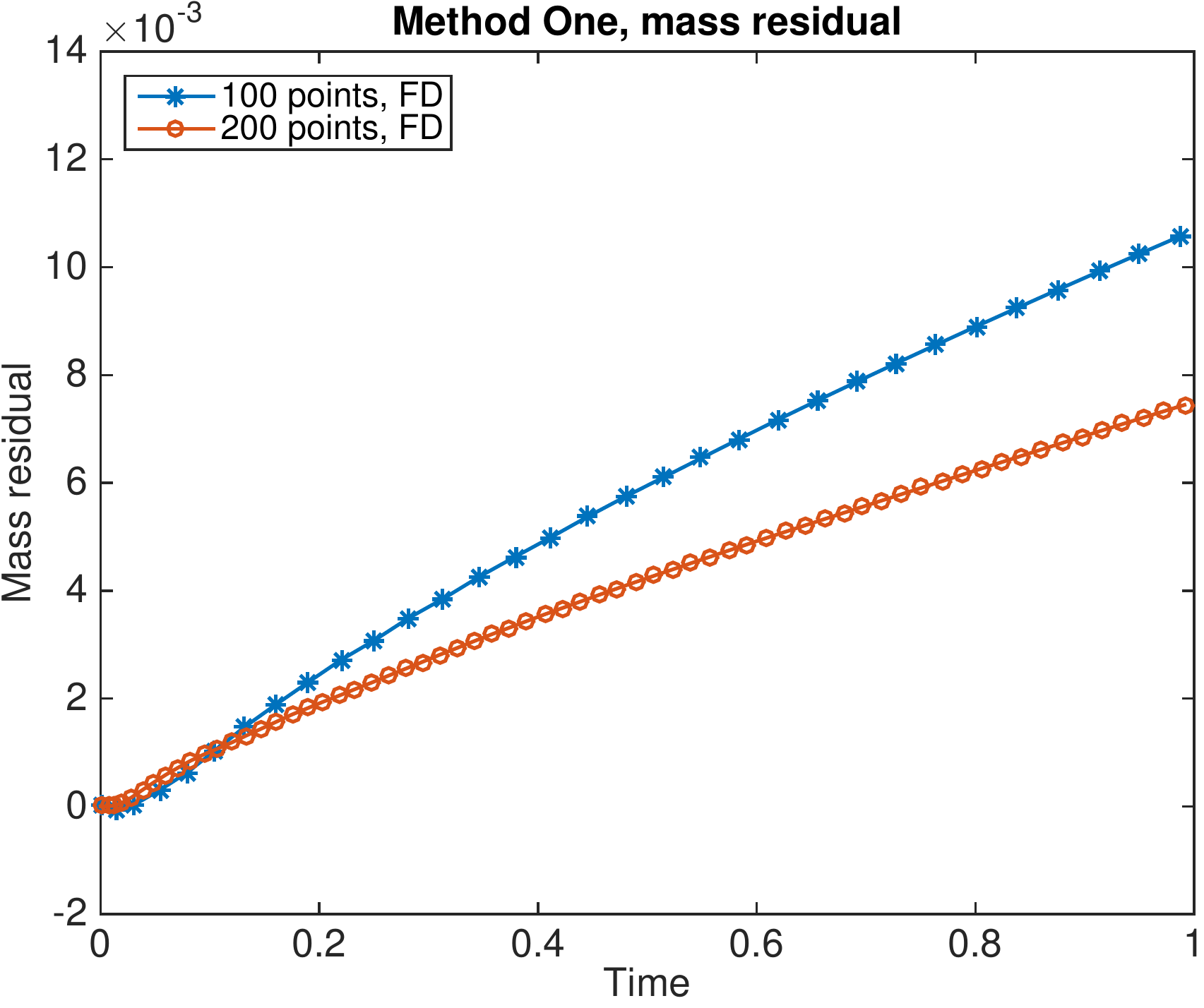} 
\includegraphics[width=0.47\textwidth ]{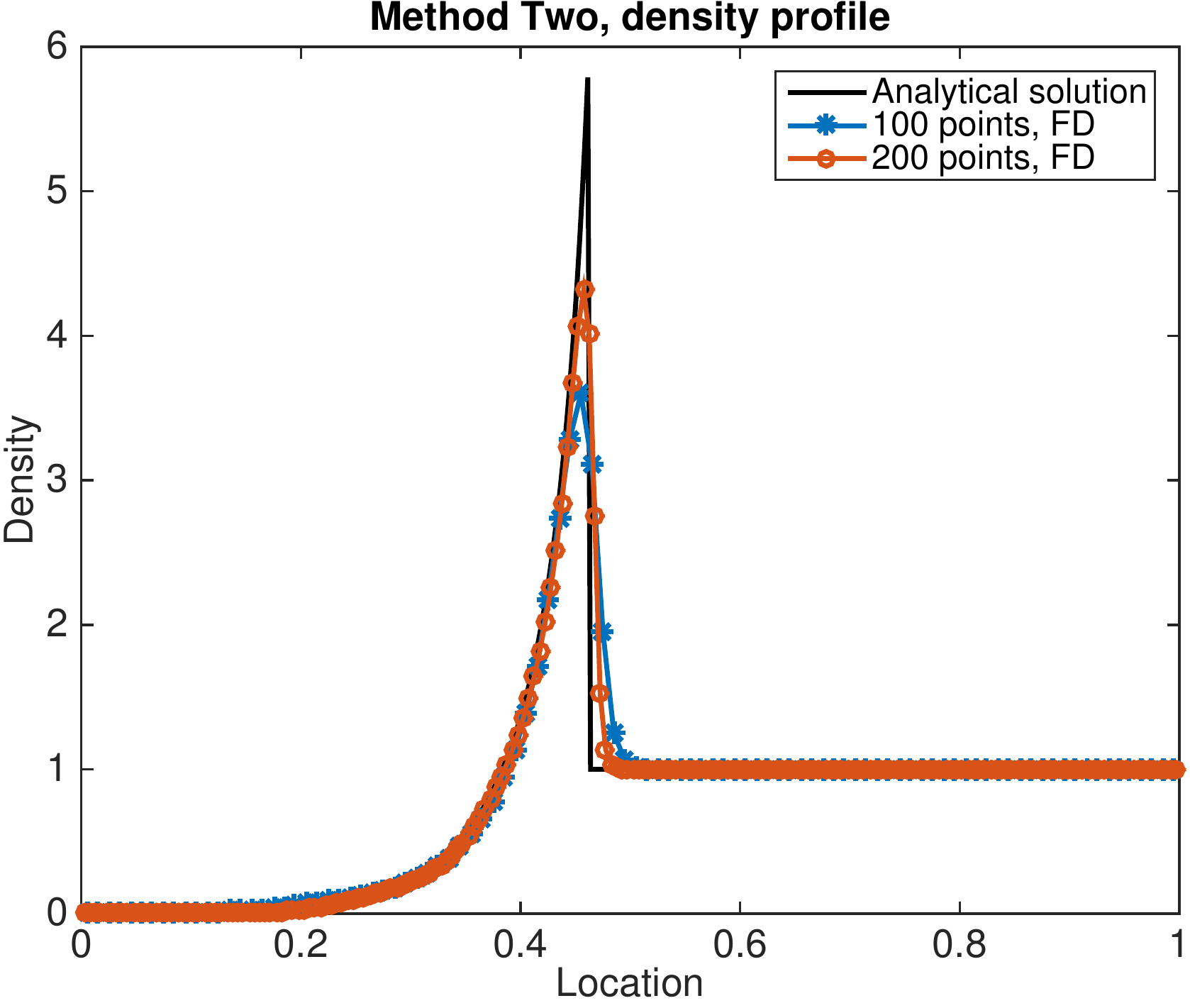} \includegraphics[width=0.48\textwidth ]{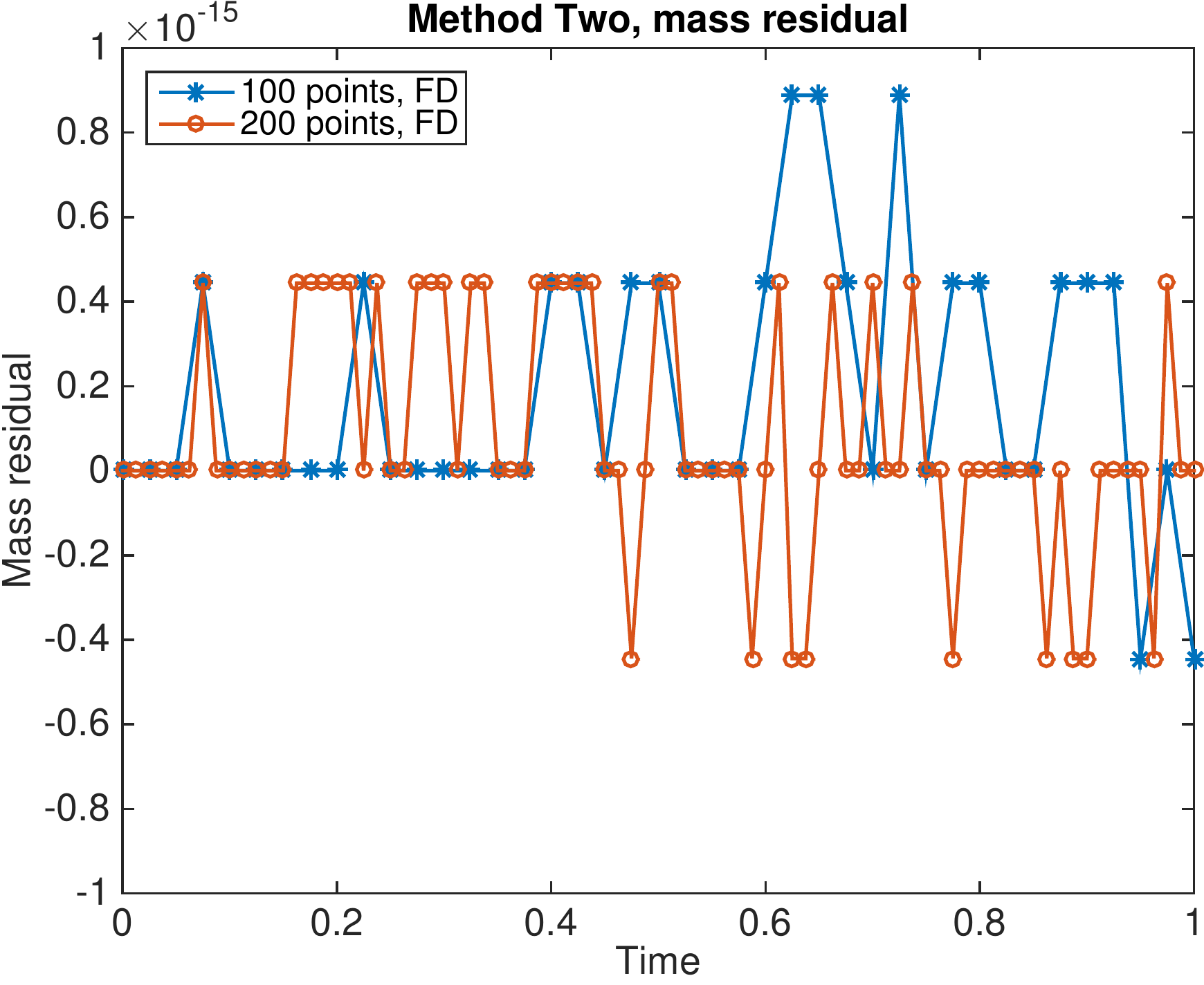}
\includegraphics[width=0.47\textwidth ]{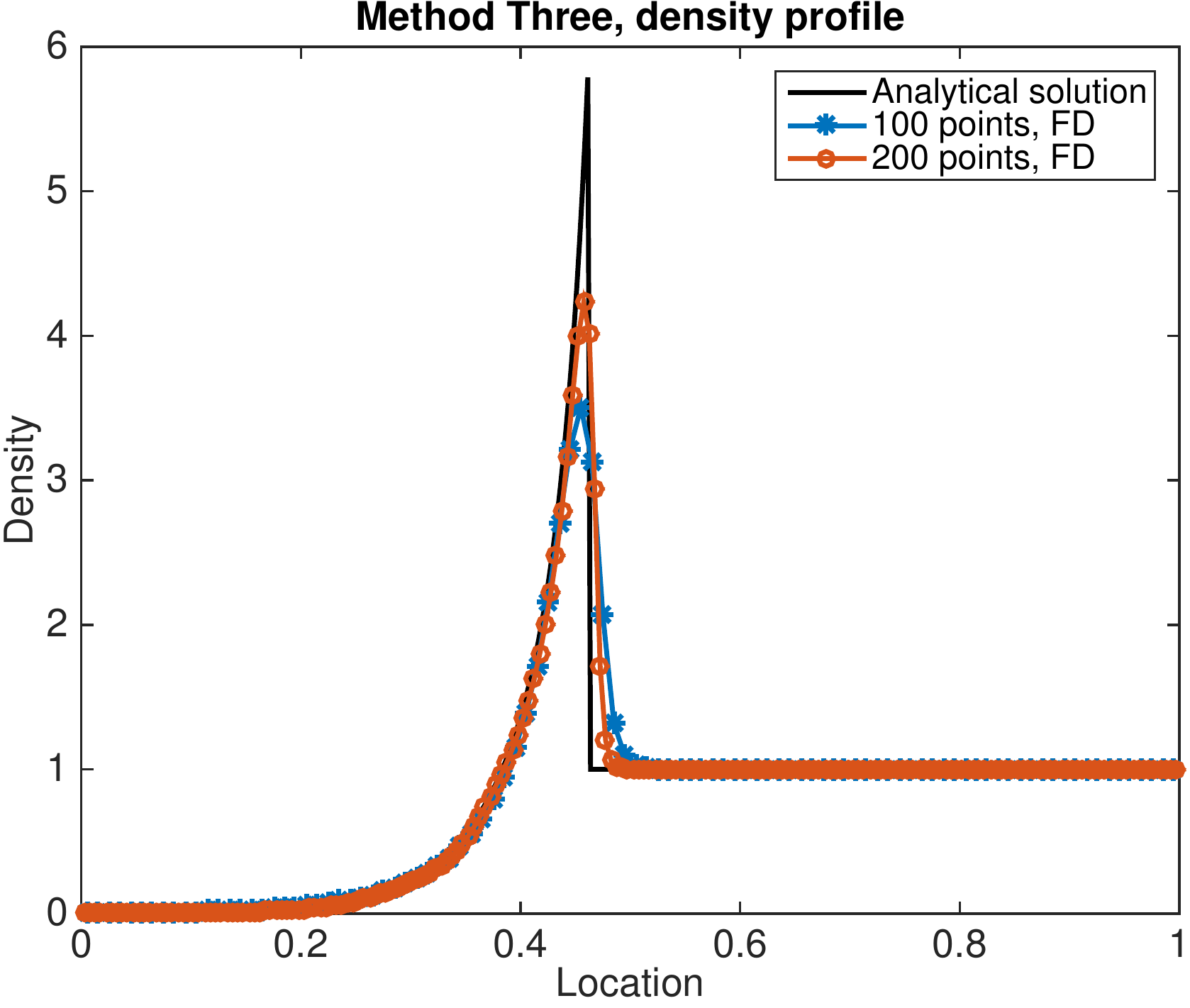} \includegraphics[width=0.48\textwidth ]{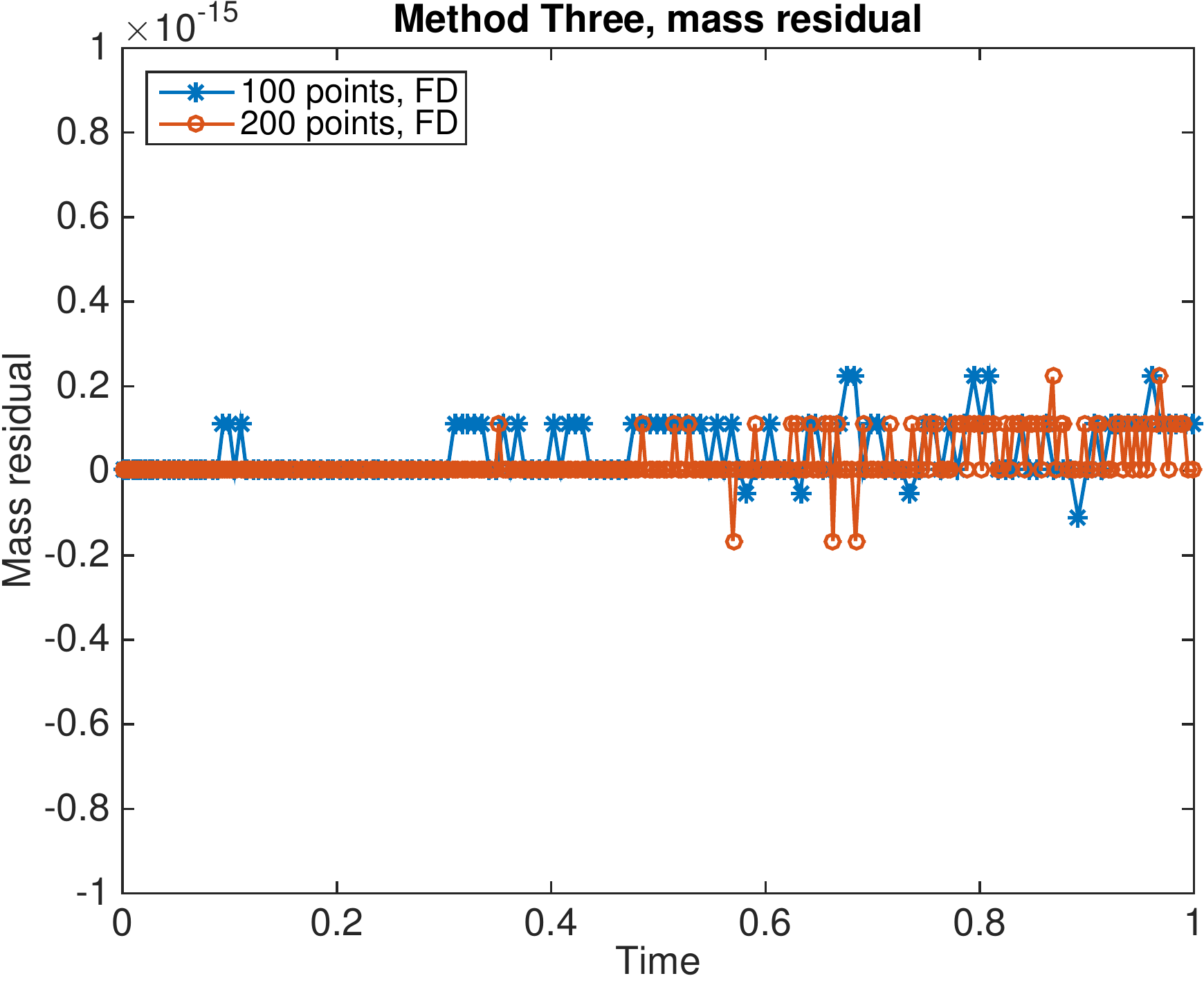}
\caption{Density profiles and mass residuals for the Sedov problem on
  different grids at t = 1.}
\label{fig:Residual of mass Sedov}
\end{figure}

   Density, velocity, and pressure profiles for  FD WENO and the analytical solution are shown in Fig.~\ref{fig:FD_Comparing Sedov}. Because the total energy residual shows results qualitatively similar to the total
   mass residual, only the latter is shown. The density profile and
   mass residual for different grid size are plotted in
   Fig.~\ref{fig:Residual of mass Sedov}.  The difference in shock
   location is clear for this problem. Method One is non-conservative
   and thus produces an incorrect shock speed and thus location; it
   appears to converge to the correct location with grid
   refinement. This result is confirmed by considering the mass
   residual. Method Two and Three can capture the right shock position on coarse gird, whereas they need much finer grid to capture the peak   of the shock. For this problem, Method Three proved to be unstable at the present Courant number due to the stiff source term, so a smaller value (0.1) is used for
   this problem.

\section{Conclusions}

We analyzed three different spatial discretizations in
cylindrical/spherical coordinates with radial dependence only for the
Euler equations using finite difference WENO. In particular, high-order accuracy and conservation were
evaluated. Only our newly proposed Method Three achieved high-order
accuracy and was conservative. The other methods are either
conservative or high-order accurate, but not both.  Current work is
underway to extend the analysis and implementations to discontinuous
finite element methods and to incorporate diffusive
effects. High-order reflecting boundary conditions will be
investigated subsequently, which are not trivial for finite
difference/volume schemes.  This approach will form the basis for
simulations of cavitation-bubble dynamics and collapse in the context
of cavitation erosion.

\section{Appendix}

\subsection*{Time discretization}
The classical fourth order explicit Runge-Kutta method for solving 
  \begin{equation}  
u_t = L(u,t)   
  \end{equation} 
  where is L(u, t) is a spatial discretization operator 
\begin{subequations}\label{eq:time discretization}
    \begin{align} 
    u^{(1)} &= u^{n} +\frac{1}{2} \Delta t L\left ( u^{(n)} \right )  \\  
    u^{(2)} &= u^{n} +\frac{1}{2} \Delta t L\left ( u^{(1)} \right )  \\ 
    u^{(3)} &= u^{n} + \Delta t L\left ( u^{(2)} \right )   \\ 
    u^{(n+1)} &=u^n+\frac{1}{6}\Delta t \left [ L(u^n) + 2 L(u^{(1)})+  2L(u^{(2)})+L(u^{(3)}) \right ] \\  \notag
    \end{align}
\end{subequations}  

All the numerical examples presented in this paper are obtained with this Runge-Kutta time discretization.

\subsection*{Local Lax-Friedrich}
The flux splitting method used in this paper is the Local Lax-Friedrich splitting scheme.

    \begin{align} \label{eq:Local Lax-Frederichs}
     f^{\pm}(u_j)=\frac{1}{2}\left ( f(u_j)\pm \lambda_j u_j  \right ) 
    \end{align}

where  $\lambda_j$ is calculated by $\lambda_j = \max\left | f^{'}(u_j) \right | $. Contrary to global Lax-Frederich, in which the wind speed is the maximum maximum absolute eigenvalue, the local Lax-Frederich use the maximum absolute eigenvalue at each point as the wind speed.

\subsection*{Implementation for the finite difference WENO scheme}
A fifth order accurate finite difference WENO scheme is applied in this paper.  For more details, we refer to \cite{Shu1996}.\\

For a scalar hyperbolic equation in 1D Cartesian coordinate
     \begin{align} 
     u_t + f(u)_x = 0  
    \end{align}

We first consider a positive wind direction, $f'(u) \geq 0$. For simplicity, we assume uniform mesh size $\Delta x$.  A finite difference spatial discretization to approximate the derivative $f(u)_x$ by

    \begin{align} 
   f(u)_x|_{x=x_{j}}\approx \frac{1}{\Delta x}\left ( \hat{f}_{j+1/2}-\hat{f}_{j-1/2}  \right ) 
    \end{align}
The numerical flux $\hat{f}_{j+1/2}$ is computed through the neighboring point values $f_j = f(u_j)$. For a 5th order WENO scheme, compute 3 numerical fluxes. The three third order accurate numerical fluxes are given by 
\begin{subequations}
  \begin{align}  
      \hat{f}_{j+1/2}^{(1)} &=\frac{1}{3}f\left ( u_{j-2} \right ) - \frac{7}{6}f\left ( u_{j-1} \right )+ \frac{11}{6}f\left ( u_{j} \right ) \\ 
      \hat{f}_{j+1/2}^{(2)} &=-\frac{1}{6}f\left ( u_{j-1} \right ) + \frac{5}{6}f\left ( u_{j} \right )+ \frac{1}{3}f\left ( u_{j+1} \right ) \\ 
      \hat{f}_{j+1/2}^{(3)} &=\frac{1}{3}f\left ( u_{j} \right ) + \frac{5}{6}f\left ( u_{j+1} \right )- \frac{1}{6}f\left ( u_{j+2} \right ) \\ \notag
 \end{align}
 \end{subequations}  
  
  The 5th order WENO flux is a convex combination of all these 3 numerical fluxes
  
  \begin{align}  
  \hat{f}_{j+1/2}= w_1 \hat{f}_{j+1/2}^{(1)} + w_2 \hat{f}_{j+1/2}^{(2)} + w_3 \hat{f}_{j+1/2}^{(3)} 
  \end{align}

and the nonlinear weights $\omega_i$ are given by
 \begin{align} 
\omega_i = \frac{\tilde{\omega}_i }{ \sum_{3}^{k=1}\tilde{\omega}_k  }, \tilde{\omega}_k = \frac{\gamma_k}{  \left ( \epsilon+\beta_k \right )} \\ 
   \end{align}
   
  with the linear weights $\gamma_k$ given by  
   \begin{align}  \notag
  \gamma_1 = \frac{1}{10},\gamma_2 = \frac{3}{5},\gamma_3 = \frac{3}{10}   \notag
     \end{align}
     
     and the smoothness indicators $\beta_k$ given by 
     \begin{subequations}
     \begin{align}  
     \beta_1 &= \frac{13}{12} \left (  f(u_{j-2})-2 f(u_{j-1})  +f(u_{j}) \right )^2 +  \frac{1}{4} \left (  f(u_{j-2})-4 f(u_{j-1})  +3f(u_{j}) \right )^2 \\
     \beta_2 &= \frac{13}{12} \left (  f(u_{j-1})-2 f(u_{j})  +f(u_{j+1}) \right )^2 +  \frac{1}{4} \left (  f(u_{j-1})-f(u_{j+1}) \right )^2  \\ 
     \beta_3 &= \frac{13}{12} \left (  f(u_{j})-2 f(u_{j+1})  +f(u_{j+2}) \right )^2 +  \frac{1}{4} \left ( 3 f(u_{j})-4f(u_{j+1}) +f(u_{j+2})\right )^2 \\\notag
     \end{align}  
     \end{subequations}    
          where $\epsilon$ is a parameter to avoid the denominator to become 0 and is usually takes as $\epsilon = 10^{-6}$ in the computation. The procedure for the case with $f'(u)$ is mirror symmetric with respect to $i+ \frac{1}{2}$.

\section*{Acknowledgments}
This work was supported in part by ONR grant N00014-12-1-0751 under Dr.\ Ki-Han Kim.

\end{document}